\documentclass[reqno, 10pt]{amsart}
\usepackage{tikz-cd}
\usepackage{enumitem}
\usepackage{amscd}
\usepackage{mathtools}
\usepackage{amssymb}
\usepackage{MnSymbol}
\usepackage[all]{xy}
\usepackage{fge}
\usepackage{hyperref}
\usepackage{graphics}
\newcommand{\moins}{\mathbin{\fgebackslash}}

\usepackage[textsize=tiny]{todonotes}

\RequirePackage[T1]{fontenc}
\RequirePackage{amsfonts,latexsym,amssymb}

\let\cal\mathcal

\usepackage{bbm}
\usepackage{hyperref}

\newtheorem{theorem}[equation]{Theorem}
 \newtheorem{lemma}[equation]{Lemma}
 \newtheorem{proposition}[equation]{Proposition}
 \newtheorem{corollary}[equation]{Corollary}  
\newtheorem{conjecture}[equation]{Conjecture}

\theoremstyle{definition}

\newtheorem{remark}[equation]{Remark}

\newtheorem{hypothesis}[equation]{Hypothesis}

\theoremstyle{remark}
\newtheorem*{acknowledgments}{Acknowledgments}

\def\invlim{\mathop{\vtop{\ialign{##\crcr$\hfill{\lim}\hfil$\crcr
\noalign{\kern1pt\nointerlineskip}\leftarrowfill\crcr\noalign
{\kern -3pt}}}}\limits}
\def\dirlim{\mathop{\vtop{\ialign{##\crcr$\hfill{\lim}\hfil$\crcr
\noalign{\kern1pt\nointerlineskip}\rightarrowfill\crcr\noalign
{\kern -3pt}}}}\limits}

\def\phi{\varphi} 
\def\epsilon{\varepsilon}

\newcommand{\ovk}{\overline{K} }

\newcommand{\dr}{\operatorname{dR} }

 \newcommand{\rank}{\operatorname{rank} }

 \newcommand{\proeet}{\operatorname{pro\acute{e}t} } 
 \newcommand{\eet}{\operatorname{\acute{e}t} }

\newcommand{\st}{\operatorname{st} }
 \newcommand{\kker}{\operatorname{Ker} } 
 \newcommand{\coker}{\operatorname{coker} }  
  \newcommand{\ccoker}{\operatorname{Coker} }  
 \newcommand{\crr}{\operatorname{cr} }
  \newcommand{\hk}{\operatorname{HK} }
   \newcommand{\LL}{\operatorname{L} }
     \newcommand{\FF}{\operatorname{FF} }
 \newcommand{\im}{\operatorname{Im} }

 \newcommand{\sv}{{\mathcal{V}}}

 \newcommand{\sll}{{\mathcal{L}}}
 
 \newcommand{\so}{{\mathcal O}}

 \newcommand{\sx}{{\mathcal{X}}}

\newcommand{\sd}{{\mathcal{D}}}

 \newcommand{\wh}{\widehat}
   \numberwithin{equation}{section}

\def\R{{\mathrm R}}
\def\O{{\cal O}}
  \def\B{{\bf B}}
\def\Q{{\bf Q}} \def\Z{{\bf Z}}

\def\O{{\cal O}}

\def\matrice#1#2#3#4{{\big(\begin{smallmatrix}#1&#2\\ #3&#4\end{smallmatrix}\big)}}
\def\lra{\longrightarrow}

\def\rg{{\rm R}\Gamma}

\def\epsilon{\varepsilon}
\def\dual{{\boldsymbol *}}
\def\bmu{{\boldsymbol\mu}}
   \newcommand{\F}{\mathbf F}

\newcommand{\Sp}{\operatorname{Sp}}

\newcommand{\dd}{\mathrm d}
\newcommand{\Xo}{X^{\circ}}
\newcommand{\So}{S^{\circ}}
\newcommand{\simto}{\xrightarrow{\;\sim\;}}

\numberwithin{equation}{section}
\begin{document}
\title[Period maps  and the $C_{\rm st}$ conjecture]
{Period maps and  the $C_{\rm st}$ conjecture}
 \author{Pierre Colmez} 
\address{CNRS, IMJ-PRG, Sorbonne Universit\'e, 4 place Jussieu, 75005 Paris, France}
\email{pierre.colmez@imj-prg.fr} 
\author{Wies{\l}awa Nizio{\l}}
\address{CNRS, IMJ-PRG, Sorbonne Universit\'e, 4 place Jussieu, 75005 Paris, France}
\email{wieslawa.niziol@imj-prg.fr}
\author{Xinyu Shao}
\address{IAZD, Leibniz Universit\"at Hannover, Welfengarten 1, 30167 Hannover, Deutschland}
\email{xinyu.shao@math.uni-hannover.de}
 \date{\today}
\thanks{The research of P.C. and W.N.
 was supported in part by the  Simons Collaboration on Perfection in Algebra, Geometry, and Topology. The research of X.S. was funded by the Deutsche Forschungsgemeinschaft (DFG, German Research Foundation) through the Collaborative Research Centre TRR 326 ``Geometry and Arithmetic of Uniformized Structures'', project number 444845124.  This material is based upon work supported by the National Science Foundation under Grant No. DMS-2424139, while W.N.  was in residence at the Simons Laufer Mathematical Sciences Institute in Berkeley, California, during the Fall 2026 semester. }
\maketitle
\begin{abstract}
We construct a  family of elliptic curves over an open disc whose
total space violates the $C_{\mathrm{st}}$ conjecture
in degree $1$. Hence the conjecture, known for proper and for
Stein varieties, as well as for products of such, does not extend to proper families over Stein bases.

This is an instance of a general criterion: the $C_{\mathrm{st}}$
conjecture holds for the total space of a smooth proper family over a
sufficiently small open disc if and only if its period map is constant,
i.e., all Gauss--Manin Higgs fields vanish or, equivalently, the $p$-adic
local systems $\R^if_*\mathbf{Q}_p$ are geometrically constant. This
criterion is special to curves: over a $2$-dimensional polydisc, we give a
family of $p$-adic analogues of non-K\"ahler elliptic surfaces whose period
maps and local systems are constant, but whose total space violates the
conjecture.
\end{abstract}
\maketitle
\tableofcontents
\section{Introduction}
Let $K$ be a complete discretely valued field of mixed characteristic
$(0,p)$ with perfect residue field $k$, and let $C=\widehat{\overline K}$.
Let  $X$ be a smooth, connected,
dagger variety over $K$.  Define 
in\footnote{We write $\Q_{p,\Box}$ for the category of solid $\Q_p$-modules.} $\sd(\Q_{p,\Box})$: 
\begin{align*}
 & \mathrm{DR}(X): =\R\Gamma_{\dr}(X)\otimes^{\LL_{\Box}}_K\B_{\dr}^+,
\quad \mathrm{DR}^i(X):=H^i(\mathrm{DR}(X))= H^i_{\dr}(X)\otimes^\Box_K\B_{\dr}^+, \\
  & F^{r,i}(X): =H^i(F^r\mathrm{DR}(X)),\quad 
   \mathrm{HK}^{r,i}(X): =\bigl(H^i_{\hk}(X_C)\otimes^{\Box}_{\breve{C}}\B_{\st}^{+}\bigr)^{N=0,\varphi=p^r},
\end{align*}
where $H^i_{\hk}(-)$ denotes the Hyodo--Kato cohomology 
and $\breve{C}=\mathrm{Frac}(W(\overline{k}))$. 
We have the following conjecture~\cite[Conj. 1.4]{CN5}, which describes the geometric pro-\'etale cohomology
of $X$ in terms of de Rham data:
\begin{conjecture}\label{main0}
For $r\geq i\geq0$, the  square
$$
\xymatrix@R=5mm{
 H^i_{\proeet}(X_C,\Q_p(r)) \ar[r] \ar[d]
   & \mathrm{HK}^{r,i}(X) \ar[d]^-{\iota_{\hk}} \\
 F^{r,i}(X) \ar[r]
   & \mathrm{DR}^i(X)
}
$$
is bicartesian in $\sd(\Q_{p,\Box}).$  Equivalently, we have the exact  sequence
\begin{equation}\label{berkeley3}
0\longrightarrow H^i_{\proeet}(X_C,\Q_p(r))
 \longrightarrow\hk^{r,i}(X)\oplus F^{r,i}(X)
 \longrightarrow \mathrm{DR}^i(X)\longrightarrow0.
\end{equation}
\end{conjecture}
For $i=0$, since $X$ is assumed to be connected,
this sequence becomes $0\to \Q_p(r)\to (\B_{\crr}^+)^{\varphi=p^r}\oplus
t^r\B_{\dr}^+\to\B_{\dr}^+\to 0$, which is one of the possible ways
to write the standard fundamental exact sequence.
Hence Conjecture~\ref{main0} is true for $i=0$, and we will assume $i\geq 1$ in the
rest of the paper.

In \cite{CN5}, this conjecture is proved for proper, Stein, and dagger-affinoid varieties, for analytifications of algebraic varieties, for almost proper varieties, and for products of a smooth proper variety with a smooth Stein variety; see \cite[Section 4]{Tate} for a recent survey.
 In \cite[Rem. 4.4]{Tate}, proper fibrations over Stein bases were proposed as a natural extension of the product case.
 In this paper we show that, without further hypotheses on the variation of the
Hodge filtration, this extension fails, and that over bases of dimension $\geq 2$
it can fail even when the Hodge filtration does not vary.

\vskip2mm
Assume $H^*_{\dr}(X)$ to be finite dimensional. If $j\leq r$,
let ${\cal E}^{r,j}(X)$ be the vector bundle on the Fargues--Fontaine curve $X_{\FF}$
associated to the $\varphi$-module\footnote{The twist $\{-r\}$ corresponds to dividing   the Frobenius by $p^r$.} $H^j_{\hk}(X_C)\{-r\}$ modified at $x_\infty$ (the point
defined by $t=0$) by imposing $\widehat\O_{X_{\FF},x_\infty}\otimes {\cal E}^{r,j}(X)=
F^r {\rm DR}^j(X)$.  Then we have an exact sequence~\cite[Cor.\,7.9]{CN5}
\begin{equation}\label{gen0}
\xymatrix@R=4mm@C=5mm{
0\ar[r] & H^1(X_{\FF},{\cal E}^{r,i-1}(X))\ar[d]
&& H^1(X_{\FF},{\cal E}^{r,i}(X))\ar[r] & 0\\
& H^i_{\proeet}(X_C,\Q_p(r))
 \ar[r] &\hk^{r,i}(X)\oplus F^{r,i}(X)
 \ar[r] & \mathrm{DR}^i(X)\ar[u] 
}
\end{equation}
So, in order to find a counterexample to Conjecture~\ref{main0} in degree~$i$,
it is enough to find $X$, such that $H^1(X_{\FF},{\cal E}^{r,i}(X))\neq 0$.

The idea behind the construction of such an example is very simple. 
Let $X\to S$
be proper and smooth and, if $s\in S$, let $X_s$ be the fiber over $s$ and
$j_s:X_s\to X$ be the corresponding closed immersion.  Assume $j_s^*:H^i_{\dr}(X)\to
H^i_{\dr}(X_s)$ is an isomorphism.  Then 
${\cal E}^{r,i}(X)\to {\cal E}^{r,i}(X_s)$ is a modification at $x_\infty$.
But ${\cal E}^{r,i}(X_s)$ is pure of slope $0$ since $X_s$ is proper; hence
$H^1(X_{\FF},{\cal E}^{r,i}(X))=0$ if and only if 
${\cal E}^{r,i}(X)\to {\cal E}^{r,i}(X_s)$ is an isomorphism, which is equivalent
to $j_s^*:F^a H^i_{\dr}(X)\to
F^a H^i_{\dr}(X_s)$ being an isomorphism for all $a$ (see Section~\ref{acyc}
for details). Hence one has to
find a family for which $j_s^*$ 
induces an isomorphism $H^i_{\dr}(X)\overset{\sim}{\to}
H^i_{\dr}(X_s)$, but not of the filtrations. 

This is achieved (Theorem~\ref{comal1})
by taking a smooth proper family of elliptic curves over a dagger disc 
for which the $H^1_{\dr}$ of the total space can be identified with that of a fiber, 
while its Hodge filtration is strictly smaller than that of the fiber (in fact, it is the part common to the Hodge filtrations of all the fibers,
see Remark~\ref{periods}).
 The same calculation applies over the open disc and gives a counterexample 
(Corollary~\ref{night1}) among smooth partially
 proper varieties that are proper over a Stein curve (a special case of holomorphically convex spaces).

\vskip2mm
  The above discussion makes it possible to give a local criterion for
Conjecture~\ref{main0} to hold for a family over a small enough disc. Let
$f:\sx\to U$ be a smooth proper family over a smooth rigid analytic curve, and let
$S\subset U$ be a sufficiently small dagger or open disc (see Hypothesis~\ref{disc} below).
Horizontal transport identifies the de Rham cohomologies of all the fibers $X_u$,
$u\in S(K)$, with a fixed space $V^i$, and the Hodge filtrations of the fibers become
a family of filtrations of $V^i$: this is the period map of the family. The Hodge
filtration of the total space $X=\sx\times_US$ is the part common to all these
filtrations (Remark~\ref{periods}), and $X$ satisfies Conjecture~\ref{main0} in
degree $i$ if and only if the period maps in degrees $i-1$ and $i$ are constant
(Corollary~\ref{period-crit}). Equivalently (Theorem~\ref{simpl2}), the Higgs fields
in degrees $i-1$ and $i$ vanish identically on $S$, or the $p$-adic local systems
$\R^{i-1}f_*\Q_p$ and $\R^if_*\Q_p$ are constant on $S_C$. For abelian schemes, the
criterion reduces to the vanishing of the Kodaira--Spencer map; for elliptic curves,
to the constancy of the $j$-invariant on $S$ (Corollary~\ref{abelian}).

\vskip2mm
This criterion is special to curves (see Section~\ref{highdim}). Over a polydisc of
dimension $\geq 2$, the conjecture still implies that the period maps and the local
systems are constant, but the converse fails: we give a family over a
$2$-dimensional polydisc whose fibers are $p$-adic analogues of non-K\"ahler
elliptic surfaces (the same construction over $\mathbf P^1$ gives the Hopf surface),
whose period maps and local systems are constant, but whose total space violates
Conjecture~\ref{main0} in degrees $2$ and $3$.  The obstruction takes values in global $2$-forms on the base, hence vanishes
when the base is a curve. (On the \'etale side, global $2$-forms on a polydisc $S$ appear as
$H^2_{\proeet}(S_C,\Q_p(2))\simeq\Omega^2(S_C)$, see~\cite[Th.\,1.8]{CDN1}). We do not know whether such
an example exists with $f$ projective.
 \begin{remark}
Conjecture 1.6 of~\cite{CN5}
allows one to extract de Rham 
and Hyodo--Kato cohomologies from pro-\'etale cohomology.  This conjecture
is actually a consequence of Conjecture~\ref{main0} when the latter is true
(see~\cite[Chap.\,9]{CN5}).
In cases where Conjecture~\ref{main0} fails, one cannot hope to recover
de Rham 
and Hyodo--Kato cohomologies via the simple recipe given in~\cite[Conj.\,1.6]{CN5},
but one must instead use the analogue of the derived formulation of~\cite{CGN} (used to recover
de Rham
and Hyodo--Kato cohomologies from {\it compactly supported} pro-\'etale cohomology).
The methods developed in~\cite{CGLN} allow one to prove the resulting statement
for quasi-compact dagger varieties, but there are subtleties with extensions
of condensed Fr\'echet spaces that prevent one from proving the 
result for (non-proper) partially proper spaces.
\end{remark}
   \begin{acknowledgments} The counterexample of Section~\ref{family} and the example of Section~\ref{highdim}
were found by GPT-6 (OpenAI).
GPT-6, Claude Fable~5.2 and Claude Opus~5.5 (Anthropic) assisted us in checking the
computations and in formulating and writing up the remaining arguments, in particular
the local criterion of Section~\ref{higgs}.  
The second author thanks GPT-6, Claude Fable~5.2 and Claude Opus~5.5 for a lot of fun!
\end{acknowledgments}

    \subsubsection*{Notation and conventions.}  Let $p$ be a fixed prime, and $\Q_p$ the field of $p$-adic numbers. 
 Let $K$ be a complete discretely valued field of mixed characteristic $(0,p)$. Let $\ovk$ be an algebraic closure of $K$ and let $C=\wh{\ovk}$ be the $p$-adic completion of $\ovk$.  Let $G_K$ be the Galois group of $K$. 
We will denote by $ \B_{\crr}, \B_{\st},\B_{\dr}$ the crystalline, semistable, and  de Rham period rings of Fontaine.

All rigid analytic spaces considered will be over $K$ or $C$. We assume that they are separated,
taut, and countable at infinity.

\section{A criterion for the conjecture to be true} 
The following criterion was used in various forms in~\cite{CN5}, 
but was not spelled out.
\begin{proposition}\label{gen1}
Let $f:X\to U$ be a morphism of dagger varieties defined over $K$. Suppose:

\quad $\bullet$ $U$ satisfies Conjecture~\ref{main0}.

\quad $\bullet$ $H^*_{\dr}(U)$ is finite dimensional
and $f^\dual:H^*_{\dr}(U)\to H^*_{\dr}(X)$ is an isomorphism.

Then $X$ satisfies Conjecture~\ref{main0}.
\end{proposition}
\begin{proof}
If ${\cal E}^{r,j}(Y)$, for $Y=U,X$, is the vector bundle associated
to $(H^j_{\hk}(Y_C)\{-r\},F^r({\rm DR}^j(Y)))$ as in the introduction, then Conjecture~\ref{main0}
for $Y$ and $i\leq r$ is equivalent to $H^1(X_{\rm FF},{\cal E}^{r,j}(Y))=0$
for $j=i-1,i$ (see~\eqref{gen0}). Now $f^*$ induces
an exact sequence $0\to {\cal E}^{r,j}(U)\to {\cal E}^{r,j}(X)\to Z^{r,j}\to 0$,
where $Z^{r,j}$ is a coherent sheaf supported at the point $x_\infty$,
which implies $H^1(X_{\FF},Z^{r,j})=0$. If $U$ satisfies Conjecture~\ref{main0},
we have 
$H^1(X_{\rm FF},{\cal E}^{r,j}(U))=0$, which implies
that $H^1(X_{\rm FF},{\cal E}^{r,j}(X))=0$, as wanted.
\end{proof}
\begin{remark}\label{simpl20}
(i) Proposition~\ref{gen1} is in fact true without the hypothesis that
$H^*_{\dr}$ is finite dimensional. To see this, one has to go back to
the basic comparison diagram~\cite[Sec.\,7.2.1]{CN5}.

(ii) As we shall see, knowing that the conjecture
is true for $X$ and that $f^\dual:H^*_{\dr}(U)\to H^*_{\dr}(X)$ is an isomorphism
 is not enough to deduce that the conjecture is true for $U$.  In fact,
if $X$ is proper, ${\cal E}^{r,j}(X)$ is of slope $0$ and, in order for $U$
to verify Conjecture~\ref{main0} for $i$, we need to have
${\cal E}^{r,j}(U)$ also of slope $0$ (otherwise some of its
slopes would be~$<0$ and $H^1(X_{\FF},{\cal E}^{r,j}(U))$
would be $\neq 0$), which is the case if and only if $f^*$
induces an isomorphism of filtered modules $H^j_{\dr}(U)\simeq H^j_{\dr}(X)$,
for $j=i,i-1$.
\end{remark}

\section{A counterexample} \label{negative1}
Let  $X$ be a smooth
dagger variety over $K$.
We define  the naive  filtration on  $\mathrm{DR}^i(X)$ as
\begin{equation}\label{ducky1}
 F^r\mathrm{DR}^i(X):=
 \im\bigl(H^i(F^r\mathrm{DR}(X))\longrightarrow H^i(\mathrm{DR}(X))\bigr). 
\end{equation}
If $H^i_{\dr}(X)$ is finite dimensional and $i\leq r$, one has also (see Proposition~\ref{appendice})
\begin{equation}\label{fil}
F^r\mathrm{DR}^i(X)=\sum_{a+b=r}F^bH^i_{\dr}(X)\otimes t^a\B_{\dr}^+
\end{equation}
If \eqref{berkeley3} is exact, the  map
\begin{equation}\label{pi}
 \pi^{r,i}: \hk^{r,i}(X)\stackrel{\iota_{\hk}}{ \longrightarrow} 
\mathrm{DR}^i(X)/F^r\mathrm{DR}^i(X)
\end{equation}
must be surjective.
That is, non-surjectivity of \eqref{pi} is enough to disprove
Conjecture~\ref{main0}.  
\begin{remark}\label{simpl21}
The relation with the proof of Proposition~\ref{gen1} is that
${\rm Coker}(\pi^{r,i})=H^1(X_{\rm FF},{\cal E}^{r,i}(X))$.
\end{remark}
\subsection{The counterexample}\label{family}
Assume $p\geq 5$, and take $K=\Q_p$. Choose $a\in \Z_p$ such that the  elliptic curve
\begin{equation}\label{fiber}
 E:\quad y^2z=x(x-z)(x-az)
\end{equation}
has good ordinary reduction (in particular, $a\not\equiv 0,1$ mod~$p$; 
any value of $a$ which does not give a CM
elliptic curve will do for a density $1$ of $p$'s; for an explicit example,
one can take $a=-3$ and $p=5$).

Let $ S=\Sp\bigl(\Q_p\langle T\rangle^\dagger\bigr), \lambda=a+p^2T, $
and let $f:X\to S$ be the projective cubic
\begin{equation}\label{family2}
 y^2z=x(x-z)(x-\lambda z),
\end{equation}
so that the fiber at $T=0$ is $j:E\hookrightarrow X$.
The dagger structure in \eqref{family2} is obtained by restriction
from the same analytic family on a disc $S'$ larger than the closed unit
disc, for example from $v_p(T)>-1$.  Since $v_p(p^2T)>0$ on that disc, both
$\lambda$ and $\lambda-1$ are units in $\O(S')$.

\begin{lemma}
The morphism $f$ is smooth and proper of relative dimension $1$, with
geometrically connected fibers and a $0$-section.  Its total space
$X$ is a smooth quasi-compact separated dagger surface.  The natural
model over $\Z_p\langle T\rangle$ has special fiber
$\overline E\times\mathbf A^1_{\F_p}$.
\end{lemma}
\begin{proof}
The discriminant of the Legendre equation is
$16\lambda^2(1-\lambda)^2$, a unit on the indicated neighborhood.
The point $[0:1:0]$ gives the section. Our lemma follows.
\end{proof}

\begin{theorem}[Counterexample]\label{comal1}
   $X$  contradicts
Conjecture~\ref{main0} in degree $1$.
\end{theorem}
\begin{proof}
It suffices  to show that the map \eqref{pi}  is not surjective for  $X$ and \(r=i=1\).
By Proposition \ref{berkeley1} below,  
$H^1_{\hk}(X_C)$ is  a rank-$2$ isocrystal over $\breve{C}$ with slopes $0,1$
and $N=0$, and $F^1\mathrm{DR}^1(X)=t\mathrm{DR}^1(X)$.  Over $\breve{C}$ choose a basis $e_0,e_1$ 
of $H^1_{\hk}(X_C)$ with $\varphi(e_0)=e_0, \varphi(e_1)=pe_1.$
The Hyodo--Kato  isomorphism identifies $\mathrm{DR}^1(X)/t\mathrm{DR}^1(X)$
 with $Ce_0\oplus Ce_1$. It follows from the proof of Proposition \ref{berkeley1} that,  
in these coordinates, $\pi^{1,1}$ becomes
$$
 \pi^{1,1}:(\B^+_{\crr})^{\varphi=p}e_0\oplus
 (\B^+_{\crr})^{\varphi=1}e_1 \longrightarrow Ce_0\oplus Ce_1,\quad (b,a)\mapsto(\theta(b),a)
$$
and is clearly not surjective (since $(\B^+_{\crr})^{\varphi=1}=\Q_p$: the cokernel is actually
$(C/\Q_p)e_1$).
\end{proof}

\begin{remark}\label{comal1.1}
(i) It follows from the above computations that ${\rm Coker}\,\pi^{1,1}$
is an almost-$C$ representation of Dimension $(1,-1)$ (it is not $C/\Q_p$
as one could infer from the above, but rather $C/\Q_p(\eta)$ for
an unramified character $\eta$, because the basis $e_0,e_1$ only exists
over $\breve C$, not $\Q_p$).

(ii) If we replace the above family by a small enough disc of curves with 
potentially multiplicative
reduction (for example, choose $a\in p\Z_p\cup (1+p\Z_p)\cup \frac{1}{p}\Z_p^\dual$, $a\neq 0,1$),
the results are exactly the same.

(iii) If, instead, we take a small enough disc of curves with supersingular reduction,
then we have ${\rm Dim}({\rm Coker}\,\pi^{1,1})=(1,-2)$, the difference
being that $\pi^{1,1}$ is injective in that case, whereas its kernel is $t\Q_p e_0$
in the other cases.
\end{remark}
\begin{remark}[Negative heights]
Let $X\to S$ be the Legendre family of Section~\ref{negative1} {\rm (}or its supersingular
variant, Remark~\ref{comal1.1}(iii){\rm )} over a sufficiently small dagger disc, and let $Y=X^{\times_S m}$, $m\geq 1$. Then we have 
\[
 \operatorname{ht} H^{2m}_{\proeet}(Y_C,\Q_p(2m))=
 \begin{cases}
  1-m  & \text{in the ordinary case},\\
  1-2m & \text{in the supersingular case}.
 \end{cases}
\]
Hence, in the ordinary case, $X\times_S X$ already gives an example of
height $-1$ in degree $4$. In the supersingular case, $X$ itself
gives an example of height $-1$ in degree $2$.

 Indeed, write $\mathcal H^j$ for the relative de Rham cohomology of
$X/S$. By relative K\"unneth, the relative de Rham cohomology of $Y/S$ in degree $2m$
is $(\mathcal H^2)^{\otimes m}$, and in degree $2m-1$ it is the direct
sum of $m$ copies of $\mathcal H^1\otimes(\mathcal H^2)^{\otimes(m-1)}$.
The class  $\ell$ of the zero  section generates $\mathcal H^2$ and is a
horizontal Tate class of Hodge degree $1$:  we have $F^1 \mathcal H^2=\mathcal H^2, F^2\mathcal H^2=0$
and $\phi(\ell)=p\ell$ in Hyodo--Kato cohomology.  Thus, by
Lemma~\ref{simpl9} below, each factor $\mathcal H^2$ contributes a Tate
twist to the absolute Hodge filtration and Hyodo--Kato cohomology. Hence, on the geometric
Fargues--Fontaine curve we  have
\[
 \mathcal E^{2m}(Y)\simeq\so_{X_{\FF}},\qquad
 \mathcal E^{2m-1}(Y)\simeq\mathcal E^1(X)^{\oplus m}.
\]
By the height formula in the proof of \cite[Prop.~7.10]{CN5}, the
desired height is the sum of the heights of
$H^0(X_{\FF},\mathcal E^{2m}(Y))$ and $H^1(X_{\FF},\mathcal E^{2m-1}(Y))$.
The first has height $1$. By Remark~\ref{simpl21}, the second is $m$ times
the height of $\ccoker\pi^{1,1}(X)$, which is $-1$ in the ordinary
case and $-2$ in the supersingular case by Remark~\ref{comal1.1}.
\end{remark}

\begin{remark}\label{comal1.2}
Consider the whole algebraic Legendre family
\[
 f\colon\mathcal E\longrightarrow
 U=\mathbb P^1_{\Q_p}\setminus\{0,1,\infty\},
 \qquad
 y^2z=x(x-z)(x-\lambda z).
\]
The total space \(\mathcal E_C^{\mathrm{an}}\)
satisfies Conjecture~\ref{main0} by
\cite[Th.~8.1]{CN5}, since it is the analytification of a
smooth algebraic variety.
Thus the conjecture holds for the whole algebraic family and fails
for its restrictions to sufficiently small dagger or open discs
by Corollary~\ref{abelian} (or Theorem~\ref{comal1} and Remark~\ref{comal1.1}).

It follows that the validity of Conjecture \ref{main0} for a variety $X$
does not imply its validity for open analytic subspaces of $X$.
\end{remark}

\begin{proposition}\label{berkeley1}
Set  $M=H^1_{\hk}(X_C)$ and $\mathrm{DR}^1={\rm DR}^1(X)$. We have
\begin{align}
  & H^1_{\dr}(X)\simto H^1_{\dr}(E),\quad  F^1H^1_{\dr}(X)=0;\label{first}\\
   & \mathrm{DR}^1\simeq(\B^+_{\dr})^2,\quad  F^1\mathrm{DR}^1=t\mathrm{DR}^1;\label{first1}\\
 & H^1_{\hk}(X_C)   \simto H^1_{\hk}(E_C),
  \quad N=0,\quad\operatorname{Slopes}_{\varphi}(M)=\{0,1\}.
 \label{slopes}
\end{align}
\end{proposition}
\begin{proof}
The first part of \eqref{first}
is proved in Proposition~\ref{unfiltered}, the second part in Lemma~\ref{ducky4};
the first part of \eqref{first1} is a consequence of Proposition~\ref{unfiltered},
the second part is Corollary~\ref{fil2}; finally \eqref{slopes} is
proved in Proposition~\ref{HK5}.
\end{proof}
\subsection{Trivialization of the Gauss--Manin connection}\label{GM}  
Let $
 {\mathcal H}:=\R^1f_*\Omega^\bullet_{X/S} $ and let $\nabla: \mathcal{H}\to \mathcal{H}\otimes_{\so_S}\Omega^1_{S/\Q_p}$  be the Gauss--Manin connection. 
We will show that  \((\mathcal H,\nabla)\)  is analytically trivial.

Let $ \omega=\dd x/(2y),\eta=x\,\dd x/(2y)$; the latter   represents a differential of the second kind.
 The classes $(\omega,\eta)$ form a basis
of ${\mathcal H}$, and the relative Hodge line is $\so_{S}\omega$.
\begin{lemma}
We have 
\begin{equation}\label{GM-matrix}
\nabla\omega
= \Big( -\frac{1}{2(\lambda-1)}\,\omega +\frac{1}{2\lambda(\lambda-1)}\,\eta \Big)\otimes d\lambda,
\quad
\nabla\eta
= \Big( -\frac{1}{2(\lambda-1)}\,\omega +\frac{1}{2(\lambda-1)}\,\eta \Big)\otimes d\lambda.
\end{equation}
Hence, the connection matrix is 
$$
 A(\lambda)=\frac{1}{2(\lambda-1)}
 \begin{pmatrix}-1&-1\\ \lambda^{-1}&1\end{pmatrix}.
$$
\end{lemma}
\begin{proof}
Work over $\Q_p(\lambda)$ and set 
$
 P(x):=x(x-1)(x-\lambda). $
We have 
\[
 \partial_\lambda\left(\frac{x^k\dd x}{2y}\right)
       =\frac{x^kx(x-1)}{4y^3}\,\dd x,\qquad k=0,1.
\]
For any rational function $R(x)$, we have
\[
 \dd_x(R/y)
 =\frac{R'P-\tfrac12RP'}{y^3}\,\dd x.
\]
Using the above, we compute: 
\begin{align*}
 \partial_\lambda\omega
 &=-\frac{1}{2(\lambda-1)}\omega
    +\frac{1}{2\lambda(\lambda-1)}\eta
    +\dd_x\left(\frac{x-x^2}{2\lambda(\lambda-1)y}\right),\\
 \partial_\lambda\eta
 &=-\frac{1}{2(\lambda-1)}\omega
    +\frac{1}{2(\lambda-1)}\eta
    +\dd_x\left(\frac{x-x^2}{2(\lambda-1)y}\right).
\end{align*}
Passing to de Rham cohomology classes yields \eqref{GM-matrix}, as wanted. 
\end{proof}

 By \eqref{GM-matrix}, for a general section $s=(\omega,\eta)w$, with $w={}^t(u,v)$ we have:
\[
\nabla(u\omega+v\eta) =
\omega\otimes du+\eta\otimes dv +u\nabla\omega+v\nabla\eta.
\]
 The section \(s\) is horizontal precisely when $
dw+A(\lambda)w\,d\lambda=0,$ that is, ${dw}/{d\lambda}=-A(\lambda)w$.
It follows that, if $G=\matrice{\alpha}{\beta}{\gamma}{\delta}$ where
$\alpha\omega+\gamma\eta, \beta\omega+\delta\eta$ form a basis of horizontal
sections, and if
\(U=\lambda-a=p^2T\), 
then $G$ satisfies 
$ G'(U)=-A(a+U)G(U)$.

Write $ A(a+U)=\sum_{n\geq0}A_nU^n. $
All the $A_n$ belong to $M_2(\Z_p)$, because $2$, $a+U$, and $a-1+U$
are units in $\Z_p[[U]]$.

\begin{lemma}\label{convergence}
Let \(G=\matrice{\alpha}{\beta}{\gamma}{\delta}\in M_2(\Q_p[[U]])\) 
be the unique formal matrix solution of
\[
G'(U)=-A(a+U)G(U),\qquad G(0)=I_2.
\]
Then \(G\) converges for \(v_p(U)>\frac{1}{p-1}\), has determinant $1$, and 
$\alpha\omega+\gamma\eta, \beta\omega+\delta\eta$
 form a horizontal basis of \(\mathcal H\) over a strict neighborhood of \(S\).
\end{lemma}
\begin{proof}
Write $G(U)=\sum_{n\geq0}G_nU^n$.  The above differential equation gives the relation
\[
 (n+1)G_{n+1}=-\sum_{j=0}^n A_jG_{n-j}
\]
from which we deduce, by induction, that
$n!G_n\in M_2(\Z_p)$:  indeed,
\[
 (n+1)!G_{n+1} =-\sum_{j=0}^n A_j\frac{n!}{(n-j)!}\bigl((n-j)!G_{n-j}\bigr).
\]
As $v_p(n!)/n\to \frac{1}{p-1}$, the series converges on the asserted disc.
The trace of $A$ is $0$, so differentiation of the determinant gives
$(\det G)'=0$ and $\det G=1$.

The series converge for $v_p(T)>\frac{1}{p-1}-2$.
In particular, $G$ and $G^{-1}$ are analytic on a disc
strictly larger than the closed unit disc.   Multiplication by
$G(p^2T)$ identifies ${\mathcal H}$, with its connection, with the constant
bundle with fiber ${\mathcal H}_0$.
\end{proof}

   The de Rham Leray spectral sequence for the smooth proper morphism $f$ is
\begin{equation}\label{leray}
 E_2^{a,b}=H^a_{\dr}\bigl(S,
                  R^bf_*\Omega^\bullet_{X/S}\bigr)
     \Longrightarrow H^{a+b}_{\dr}(X).
\end{equation}
Here cohomology on the base is computed using the Gauss--Manin connection.  

\begin{proposition}\label{unfiltered}
Restriction to $T=0$ induces isomorphisms
\[
 j^*:H^i_{\dr}(X)\simto H^i_{\dr}(E)
 \quad\hbox{for every }i.
\]  
In particular, the dimensions
in degrees $0,1,2$ are $1,2,1$, and the other groups vanish. 
The same claim holds after base change from \(\Q_p\) to \(C\).
\end{proposition}
\begin{proof}
The relative cohomology bundles are $0$ outside degrees $0,1,2$.
In degree $0$ the bundle is just $(\so,\dd)$; in degree $1$ it is
constant as a connection by Lemma~\ref{convergence}; and in degree
$2$ it is a trivial rank-$1$ connection, generated by the relative
class of the $0$-section.  Since $H^i_{\dr}(S)=0$, for $i>0$, each of
these three connections has $H^1_{\dr}=0$ on the
base.  Thus \eqref{leray} is supported in the column $a=0$.
Evaluation of a horizontal section at $T=0$ identifies that column
with the cohomology of $E$, proving the proposition.
\end{proof}
\begin{remark} In degree $1$, 
the isomorphisms in Proposition~\ref{unfiltered} are not 
strict for the Hodge filtration; we will show this 
in the next computation.
\end{remark}
\subsection{The Hodge filtration}\label{hodge}
We will show that $F^1H^1_{\dr}(X)=0$. 
We have  an exact sequence
$$
0\longrightarrow f^*\Omega^1_{S/\Q_p}
\longrightarrow\Omega^1_{X/\Q_p}
\longrightarrow\Omega^1_{X/S}\longrightarrow0.
$$
Writing $\omega_f=f_*\Omega^1_{X/S}$ and using
$f_*\so_{X}=\so_{S}$, its pushforward starts as
\begin{equation}
\begin{split}
0\longrightarrow\Omega^1_{S/\Q_p}
\longrightarrow f_*\Omega^1_{X/\Q_p}
\longrightarrow\omega_f
\xrightarrow{\mathrm{KS}}
(R^1f_*\so_{X})\otimes\Omega^1_{S/\Q_p},
\end{split}
\end{equation}
where $\mathrm{KS}$ denotes  the Kodaira--Spencer map.

\begin{lemma}\label{lem:KS}
The  map $\mathrm{KS}$  is an isomorphism of line bundles. Hence,
$ f_*\Omega^1_{X/\Q_p}\simeq \Omega^1_{S/\Q_p}.$
\end{lemma}
\begin{proof}The Kodaira--Spencer map is the composite (we identify \(\mathcal H/\omega_f\simeq \R^1f_*\mathcal O_X\))
\[
\omega_f \to \mathcal H \xrightarrow{\nabla} \mathcal H\otimes\Omega^1_{S/\Q_p}
\xrightarrow{\mathrm{pr}\otimes1} (\mathcal H/\omega_f)\otimes\Omega^1_{S/\Q_p}.
\] 
The first formula  of
\eqref{GM-matrix} gives 
\[
 \nabla\omega\bmod\omega_f
   =\frac{1}{2\lambda(\lambda-1)}\,\overline\eta\otimes\dd\lambda
   =\frac{p^2}{2\lambda(\lambda-1)}\,\overline\eta\otimes\dd T.
\]
 Our lemma follows since  ${p^2}/({2\lambda(\lambda-1)})$ is a unit in $\so(S)$. 
\end{proof}

\begin{lemma}\label{ducky4}
We have $F^1H^1_{\dr}(X)=0$.
\end{lemma}
\begin{proof}
We have $ F^1H^1_{\dr}(X)={\rm Im}(\Omega^1(X)^{d=0}\to H^1_{\dr}(X)).$
By Lemma~\ref{lem:KS}, every global $1$-form is pulled back
from $S$.  Such a form is closed since $S$ is a curve, and is exact since $H^1_{\dr}(S)=0$;
its de Rham class
on $X$ is therefore $0$.  
\end{proof}
\begin{remark}
There is no contradiction between Lemma \ref{ducky4} and the $1$-dimensional \(F^1\) on each fiber:
a relative invariant differential does not extend to a global $1$-form on $X$
when the Kodaira--Spencer map is nonzero.
\end{remark}

\begin{corollary}\label{fil2}
$F^1{\rm DR}^1=t {\rm DR}^1$.
\end{corollary}
\begin{proof}
This is a direct consequence of~\eqref{fil}, granted Lemma~\ref{ducky4}.
\end{proof}
\subsection{Hyodo--Kato cohomology}
We  have  the functorial and compatible Hyodo--Kato isomorphisms (see \cite[Prop. 4.24, Prop. 4.27, Cor. 5.12]{CN5}) 
\begin{align}\label{HK1}
  & H^i_{\hk}(X_C)\otimes^{\Box}_{\breve{C}}C\simto H^i_{\dr}(X_C),\\
 &  H^i_{\hk}(X_C)\otimes^{\Box}_{\breve{C}}\B^+_{\dr}
                \simto H^i_{\dr}(X)\otimes^{\Box}_K \B^+_{\dr}.\notag
\end{align}

\begin{proposition}\label{HK5}
$j:E_C\hookrightarrow X_C$ induces an
isomorphism of $(\varphi,N)$-modules:
\[
 j^*:H^i_{\hk}(X_C)\simto H^i_{\hk}(E_C)
 \quad\hbox{for every }i.
\]
\end{proposition}
\begin{proof}
By functoriality of \eqref{HK1} this follows from 
Proposition~\ref{unfiltered}.  
\end{proof}
Combined with the hypothesis that $E$ has good ordinary reduction, this yields 
\eqref{slopes} in Proposition \ref{berkeley1}.

\begin{remark}[Restriction to the open disc]\label{open}
Let $\So=\{v_p(T)>0\}$  and let $\Xo\to\So$ be the 
restriction of \eqref{family2}.  This is a smooth proper elliptic
fibration over a Stein curve, and the total space is partially proper.
Since the Poincar\'e lemma holds for the open disc,
the same arguments as in the preceding proofs yield the following result.
\end{remark}

 \begin{corollary}\label{night1}
$\Xo$  contradicts Conjecture \ref{main0} in degree $1$.
\end{corollary}

\subsection{Alternative proof of Corollary \ref{night1}}
We now compute $H^1_{\proeet}(\Xo_C,\Q_p(1))$ directly, using Kummer theory and Serre--Tate theory, and derive from that another proof of Corollary \ref{night1}.

\begin{proposition}\label{prop:etale}
Let $\sv=\R^1f_{*}\Q_p(1)$, a $\Q_p$-local system on $\So_C$. 

{\rm (i)} We have $H^0(\So_C,\sv)\simeq\Q_p(1)$.

{\rm (ii)} There is an exact sequence:
\begin{equation}\label{etale-seq}
 0\longrightarrow \so(\So_C)/C
 \longrightarrow H^1_{\proeet}(\Xo_C,\Q_p(1))
 \longrightarrow \Q_p(1)
 \longrightarrow 0.
\end{equation}
\end{proposition}
\begin{proof}
\emph{Step 1: The implication {\rm (i)}$\Rightarrow${\rm (ii)}.}  
Since $f$ is proper and smooth, the sheaves
$\R^bf_{*}\Q_p(1)$ are $\Q_p$-local systems on $\So_C$ \cite[Th. 10.5.1]{SW20}.
We have $\R^0f_{*}\Q_p(1)=\Q_p(1)$, because the fibers of $f$ are geometrically
connected, and 
$\sv:=\R^1f_{*}\Q_p(1)$ 
is the rational Tate module $V_p(\Xo_C/\So_C)$ of the family.  
Consider the Leray spectral sequence
\[
 E_2^{a,b}=H^a_{\proeet}(\So_C,\R^bf_{*}\Q_p(1))
 \Longrightarrow H^{a+b}_{\proeet}(\Xo_C,\Q_p(1)).
\]
Since $H^2_{\proeet}(\So_C,\Q_p(1))=0$ (because $\So$ is a Stein curve, see~\cite[Th.\,1.8]{CDN1}), 
we get an exact sequence
\begin{equation}\label{leray-et}
 0\longrightarrow H^1_{\proeet}(\So_C,\Q_p(1))
 \xrightarrow{\ f^*\ } H^1_{\proeet}(\Xo_C,\Q_p(1))
 \longrightarrow H^0(\So_C,\sv)\longrightarrow0.
\end{equation}
Since $H^i_{\mathrm{dR}}(S_C^\circ)=0$ for $i>0$, we have $H^1_{\mathrm{pro\acute et}}(S_C^\circ,\mathbf Q_p(1))
\simeq \mathcal O(S_C^\circ)/C$
by~\cite[Th.\,1.8]{CDN1}. Thus (ii) follows from~(i) and~\eqref{leray-et}.  
Hence we just have to prove (i).

\emph{Step 2: An integral model.}
Put $W=W(\overline{\F}_p)$ and $R=W[[T]]$, with its $(p,T)$-adic
topology. The Legendre equation with $\lambda=a+p^2T$ defines an
elliptic curve $\mathfrak X/R$ with ordinary closed fiber $\overline E$.
Fix identifications of the multiplicative and \'etale parts of
$\overline E[p^\infty]$ with $\bmu_{p^\infty}$ and $\Q_p/\Z_p$.
Then $G=\mathfrak X[p^\infty]$ is the extension
\begin{equation}\label{conn-et}
 0\longrightarrow\bmu_{p^\infty}\longrightarrow G
 \longrightarrow\Q_p/\Z_p\longrightarrow 0
\end{equation}
with Serre--Tate parameter $q\in 1+(p,T)R$
\cite[Th.~2.1]{Ka81}. In particular, $q$ is an analytic unit on $\So_C$.

\emph{Step 3: The extension class.}
The rational Tate module of~\eqref{conn-et} gives
\[
 0\longrightarrow\Q_p(1)\longrightarrow\sv
 \longrightarrow\Q_p\longrightarrow 0
\]
on $\So_C$. Its boundary sends $1$ to the Kummer class $\kappa(q)$,
since the fiber of $G[p^n]\to\Z/p^n$ over $1$ is the torsor of
$p^n$-th roots of $q$. The natural map
$H^1_{\proeet}(\So_C,\Q_p(1))\to\Omega^1(\So_C)$ sends $\kappa(q)$
to $\dd q/q$~\cite{CDN1}. Hence, if $q$ is nonconstant, the boundary
is injective and $H^0(\So_C,\sv)=\Q_p(1)$.

 \emph{Step 4:  Nonconstancy of $q$.}
Since the Legendre $j$-invariant $J(\lambda)$ is a nonconstant rational
function, choose $t\in pW$ such that $J(a+p^2t)\ne J(a)$. The fibers
$\mathfrak X_t$ and $\mathfrak X_0$ are lifts of the same $\overline E$
to $W$, with Serre--Tate parameters $q(t)$ and $q(0)$, respectively.
If $q(t)=q(0)$, these lifts would be isomorphic by
\cite[Th.~2.1]{Ka81}, contradicting that their $j$-invariants are not equal.
Thus $q$ is nonconstant on $\So_C$, proving (i).
\end{proof}
\begin{remark}\label{rem:etale}
(1) By \eqref{etale-seq}, the height of\footnote{To justify
this statement, one has to upgrade \eqref{etale-seq} to the category of qBC's or define the height
for their almost-$C$-representation avatars (if $K$ is a finite extension of $\Q_p$).
Recall that, for an almost $C$-representation $W$, 
we have $\sum_{i=0}^2(-1)^i\dim_{\Q_p}H^i(G_K,W)=-[K:\Q_p]{\rm ht}(W)$, see~\cite{Fo-Cp}
(this follows from Tate's Galois cohomology computations). In our case,
this would make no sense since $H^0$ and $H^1$ are infinite dimensional.
But $\cup\log\chi_{\rm cycl}$ induces a map
$\ell:H^0(G_K,W)\to H^1(G_K,W)$, with finite dimensional kernel and cokernel.
Hence we can define ${\rm ht}(W)$ by the formula
$\dim_{\Q_p}{\rm Coker}\,\ell-\dim_{\Q_p}{\rm Ker}\,\ell-\dim_{\Q_p}H^2(G_K,W)=
[K:\Q_p]{\rm ht}(W)$.}
 $H^1_{\proeet}(\Xo_C,\Q_p(1))$
is $1$, whereas, if Conjecture~\ref{main0} were true,
 it should be $\dim_CH^1_{\dr}(\Xo_C)=2$
according to \cite[Prop.~7.10]{CN5}.  
This gives a second
proof of Corollary~\ref{night1}. 

(2) The two computations agree.  The basic comparison sequence
\cite[Cor.~7.9]{CN5}, \cite{CN4}  in degree $1$, for $r=1$, reads 
\[
 \hk^{1,0}(\Xo)\oplus F^{1,0}(\Xo)\to\mathrm{DR}^0(\Xo)\to
 H^1_{\proeet}(\Xo_C,\Q_p(1))\to\hk^{1,1}(\Xo)\oplus F^{1,1}(\Xo)\to\mathrm{DR}^1(\Xo).
\]
The first map is surjective, since $\mathrm{DR}^0(\Xo)=\B^+_{\dr}$,
$F^{1,0}(\Xo)=t\B^+_{\dr}$, and
$\theta\bigl((\B^+_{\crr})^{\varphi=p}\bigr)=C$.  The kernel of
$F^{1,1}(\Xo)\to\mathrm{DR}^1(\Xo)$ is the cokernel of
$H^0(\mathrm{DR}(\Xo))\to H^0(\mathrm{DR}(\Xo)/F^1)=\so(\Xo_C)$, i.e.,
it is $\so(\So_C)/C$; its image is $F^1\mathrm{DR}^1=t\,\mathrm{DR}^1$ by Remark \ref{open}. 
 We get the exact sequence
\[
 0\longrightarrow\so(\So_C)/C\longrightarrow H^1_{\proeet}(\Xo_C,\Q_p(1))
 \longrightarrow\kker\pi^{1,1}\longrightarrow0,
\]
and $\kker\pi^{1,1}=\Q_p\,t\,e_0\simeq\Q_p(1)$ by the proof of
Theorem~\ref{comal1}.  Thus the two computations give the same
extension.  On the \'etale side, the quotient $\Q_p(1)=H^0(\So_C,\sv)$ is
the rational Tate module of the multiplicative part $\bmu_{p^\infty}$ of $G$.
  On the de Rham side, it is $\kker\pi^{1,1}$.
\end{remark}

\section{A local criterion}\label{higgs}
We now give a necessary and sufficient condition for Conjecture~\ref{main0}
to hold for a family over a small disc.

 Let $U$ be a smooth rigid analytic curve over $K$, and let
$f:\mathcal X\to U$ be smooth and proper. By Scholze's relative
comparison theorem \cite[Th.~8.8]{Sch13}, whose hypothesis on the
higher direct images of the constant $\widehat{\Z}_p$-local system
follows from \cite[Th.~10.5.1]{SW20},
the relative Hodge--de Rham spectral sequence
\[
 E_1^{a,b}=R^bf_*\Omega^a_{\mathcal X/U}
 \Longrightarrow
 \mathcal H^{a+b}:=R^{a+b}f_*\Omega^\bullet_{\mathcal X/U}
\]
degenerates at $E_1$, and its terms are locally free. They commute
with base change by proper coherent base change. Write
$F^\bullet\mathcal H^i$ for the relative Hodge filtration and
$\nabla$ for the Gauss--Manin connection.
Griffiths transversality induces the Higgs maps
\[
 \theta_i^a:\quad \operatorname{gr}_F^a\mathcal H^i
 \longrightarrow
 \operatorname{gr}_F^{a-1}\mathcal H^i\otimes_{\so_U}\Omega^1_{U/K}.
\]
Their vanishing is equivalent to
\begin{equation}\label{higgs0}
 \nabla(F^a\mathcal H^i)
 \subseteq F^a\mathcal H^i\otimes_{\so_U}\Omega^1_{U/K},
 \qquad\text{for all }a,i.
\end{equation}

\subsection{The Hodge filtration}\label{last1}
Let $S$ be a dagger affinoid in $U$ or a Stein open subspace of $U$, and set $X=\sx\times_U S$.

\begin{lemma}\label{lift}
For every $a,i\geq0$, the  edge map induces a surjection
\begin{equation}\label{eq:filtered-lift}
 H^i(X,\Omega^{\geq a}_{X/K})
 \longrightarrow
 \Gamma(S,F^a\mathcal H^i)^{\nabla=0}.
\end{equation}
\end{lemma}
\begin{proof}
Since the base is a curve, we have an 
 exact sequence 
\begin{equation}\label{higgs-trunc}
 0\longrightarrow
 f^{-1}\Omega^1_{S/K}\otimes_{f^{-1}\so_S}
       \Omega^{\geq a-1}_{X/S}[-1]\longrightarrow \Omega^{\geq a}_{X/K}
 \longrightarrow \Omega^{\geq a}_{X/S}
 \longrightarrow0.
\end{equation}
Relative Hodge--de Rham degeneration identifies $ \R^if_*\Omega^{\geq a}_{X/S}\simto F^a\mathcal H^i
 \subseteq\mathcal H^i.$
Hence, coherent acyclicity on $S$ identifies the relevant part
of the cohomology sequence of \eqref{higgs-trunc} with 
$$
H^i(X,\Omega^{\geq a}_{X/K})
 \longrightarrow H^0(S,F^a\mathcal H^i)\xrightarrow{\nabla}
 H^0(S,F^{a-1}\mathcal H^i\otimes_{\so_S}\Omega^1_{S/K}).
$$
The  map $\nabla$  is the filtered Gauss--Manin connection.  Its target
injects into
$H^0(S,\mathcal H^i\otimes_{\so_S}\Omega^1_{S/K})$.
Our lemma follows.
\end{proof}

 \begin{hypothesis}\label{disc}
We fix $s_0\in U(K)$, and $S\subset U$ is a dagger disc or an open disc centered at $s_0$,
small enough for each of the finitely many connections $({\cal H}^i,\nabla)$ to admit
a horizontal trivialization on a neighborhood of $S$:
\begin{equation}\label{G-triv}
 {\cal H}^i|_S\simeq \O_S\otimes H^0(S,{\cal H}^i)^{\nabla=0}.
\end{equation}
We set $X:=\sx\times_US$ {\rm (}if $S$ is a dagger disc, $X$ is a dagger variety, whose
dagger structure comes from the $\sx\times_US'$, for $S'$ running through the discs
centered at $s_0$ of radius strictly bigger than that of $S${\rm )}.
\end{hypothesis}
Such discs exist: the argument is the same as for Lemma~\ref{convergence},
i.e., choose $S_0$ a small enough disc around $s_0$ such that 
${\cal H}^i$ admits a basis over $\O(S_0)$; write the equation $\nabla=0$
in this basis;  one gets a differential equation $G'=-AG$, and the formal solution
with $G(s_0)=1$ converges and is invertible on some disc
around $s_0$.  We will use the following properties of such an $S$: it is geometrically connected,
$S(K)$ is Zariski dense in $S$, coherent sheaves on $S$ are acyclic, and
$H^i_{\dr}(S)=0$ for $i>0$ (Poincar\'e lemma).

If $s\in S(K)$, denote by $X_s$ the fiber above $s$, and by $j_s$ the closed
embedding $X_s\hookrightarrow X$.

\begin{lemma}\label{simpl9}Let $S$ and $X$ be as in Hypothesis~\ref{disc}.

{\rm (i)} $F^a H^i_{\dr}(X)\simeq H^0(S,F^a{\cal H}^i)^{\nabla=0}$ for all $a$.

{\rm (ii)}
$j_s^\dual$
induces isomorphisms $j_s^\dual:H^i_{\dr}(X)\overset{\sim}{\to} H^i_{\dr}(X_s)$
and $j_s^\dual:H^i_{\hk}(X_C)\overset{\sim}{\to} H^i_{\hk}(X_{s,C})$.
\end{lemma}
\begin{proof}
Set $V^j:=H^0(S,{\cal H}^j)^{\nabla=0}$.
Consider the Leray spectral sequence
$E_2^{i,j}=H^i_{\dr}(S,\mathcal H^j)\Rightarrow H^{i+j}_{\dr}(X)$.
By \eqref{G-triv}, $H^i_{\dr}(S,\mathcal H^j)=H^i_{\dr}(S)\otimes_KV^j$.
This vanishes for $i>0$ by the Poincar\'e lemma for the dagger and open discs.
Hence the edge map
\[
 e:\ H^i_{\dr}(X)\longrightarrow H^0_{\dr}(S,\mathcal H^i)
   =H^0(S,\mathcal H^i)^{\nabla=0}
\]
is an isomorphism.  Its composite with evaluation at $s$ is $j_s^*$, which proves (i)
for $a=0$.

The map $e$ is induced by $\Omega^\bullet_{X/K}\to\Omega^\bullet_{X/S}$,
together with the identification
$H^i(X,\Omega^\bullet_{X/S})=H^0(S,\mathcal H^i)$ given by coherent
acyclicity on $S$.  Consequently, the restriction of $e$ to the image of
$H^i(X,\Omega^{\geq a}_{X/K})$ is induced by
$\Omega^{\geq a}_{X/K}\to\Omega^{\geq a}_{X/S}$.  That is, it is the edge map
\eqref{eq:filtered-lift} followed by the inclusion
$H^0(S,F^a\mathcal H^i)\subseteq H^0(S,\mathcal H^i)$.  By
Lemma~\ref{lift}, its image is
 $H^0(S,F^a\mathcal H^i)^{\nabla=0}$.
This proves (i).

For $H^i_{\dr}$, (ii) follows from (i), \eqref{G-triv} and the fact that
the fiber at $s$ of ${\cal H}^i$ is $H^i_{\dr}(X_s)$, by construction.
The statement for $H^i_{\hk}$ follows from that for $H^i_{\dr}$ and
functoriality of the Hyodo--Kato isomorphism. 
\end{proof}

\subsection{Acyclicity of vector bundles and vanishing of the Higgs map}\label{acyc}
Let $${\cal E}^{i}(X):={\cal E}(H^i_{\hk}(X_C),F^0( H^i_{\dr}(X)\otimes_K\B_{\dr}))$$
be the\footnote{The vector bundle ${\cal E}^i(X)$ is isomorphic to
the vector bundle ${\cal E}^{r,i}(X)$ defined in the introduction; they only differ
by a Tate twist if one takes Galois action into account.}
 vector bundle on the Fargues--Fontaine curve $X_{\rm FF}$ attached
to $(H^i_{\hk}(X_C),H^i_{\dr}(X))$.
Conjecture~\ref{main0} for $r,i$ is equivalent to
$H^1(X_{\FF},{\cal E}^{j}(X))=0$, for $j=i-1,i$; i.e., ${\cal E}^j(X)$ {\it is acyclic}
for $j=i,i-1$.

Now, by Lemma~\ref{simpl9}, $j_s^\dual$
induces isomorphisms $j_s^\dual:H^i_{\dr}(X)\overset{\sim}{\to} H^i_{\dr}(X_s)$
and $j_s^\dual:H^i_{\hk}(X_C)\overset{\sim}{\to} H^i_{\hk}(X_{s,C})$,
and sends $F^a H^i_{\dr}(X)$ to $F^a H^i_{\dr}(X_s)$ (but the
induced map has no reason to be an isomorphism, a priori).
This induces an exact sequence of coherent sheaves on $X_{\rm FF}$:
$$0\to {\cal E}^{i}(X)\to {\cal E}^{i}(X_s)\to Z^i\to 0$$
where $Z^i$ is supported on the point $x_\infty$
(defined by the equation $t=0$), with
$$H^0(X_{\FF},Z^i)=
F^0(H^i_{\dr}(X_s)\otimes_K\B_{\dr})/F^0(H^i_{\dr}(X)\otimes_K\B_{\dr})$$
(a $\B_{\dr}^+$-module killed by a power of $t$).

But $X_s$ is proper, hence ${\cal E}^{i}(X_s)$ is pure, of slope $0$,
so that $H^0(X_{\rm FF},{\cal E}^{i}(X_s))$ is finite dimensional over $\Q_p$,
and $H^1(X_{\rm FF},{\cal E}^{i}(X_s))=0$.
It follows that ${\cal E}^i(X)$ is acyclic if and only if
$H^0(X_{\FF},Z^i)=0$
or,
equivalently, if and only if $j_s^\dual:F^a H^i_{\dr}(X)\to F^a H^i_{\dr}(X_s)$
is an isomorphism, for all $a$.
Hence we have proved the following result:
\begin{lemma}\label{simpl1}
The following conditions are equivalent:

{\rm (i)} ${\cal E}^i(X)$ is acyclic.

{\rm (ii)} There exists $s\in S(K)$, such that $j_s^\dual:F^a H^i_{\dr}(X)\to F^a H^i_{\dr}(X_s)$
is an isomorphism, for all~$a$.

{\rm (iii)} For all $s\in S(K)$, $j_s^\dual:F^a H^i_{\dr}(X)\to F^a H^i_{\dr}(X_s)$
is an isomorphism, for all~$a$.
\end{lemma}

It follows that condition (iii) of Lemma~\ref{simpl1} for $a$ is also equivalent to:

\quad (iv) The natural injection $\O_S\otimes  H^0(S,F^a{\cal H}^i)^{\nabla=0}\to
F^a{\cal H}^i$ is an isomorphism.

(That it is injective is a general fact in the theory of differential equations;
that it is surjective if and only if (iii) is true follows from the fact that the fiber of
$F^a{\cal H}^i$ at $s$ is $F^aH^i_{\dr}(X_s)$.)

Now, since $\O_S\otimes  H^0(S,{\cal H}^i)^{\nabla=0}\to
{\cal H}^i$ is an isomorphism, (iv) is equivalent to

\quad (v) $F^a{\cal H}^i$ is stable by $\nabla$.

This in turn is equivalent to the vanishing of the 
Higgs map $\theta_i^a$.
Putting everything together, we obtain the equivalence between (i)
  and (ii) in the following theorem:
\begin{theorem}\label{simpl2}
Let $S$ and $X$ be as in Hypothesis~\ref{disc}
and let $i\geq 1$.
The following conditions are equivalent:

{\rm (i)} Conjecture~\ref{main0} is true for $X$ in degree $i$
{\rm (}for one or, equivalently, for all $r\geq i${\rm )}.

{\rm (ii)} All Higgs maps $\theta_j^a$ are $0$ on $S$, for $j=i-1,i$.

{\rm (iii)} The $\Q_p$-local systems $\R^jf_*\Q_p$ are constant
on $S_C$\footnote{If $S$ is a dagger disc, this means: on $S'_C$, for some disc $S'\subset U$ centered at $s_0$, of radius strictly bigger than that of $S$.}, for $j=i-1,i$.
\end{theorem}

\begin{proof}
The equivalence between (i) and (ii) has been proved above; we are going to prove
(i)$\Rightarrow$(iii) and (iii)$\Rightarrow$(ii).
Fix $j\in\{i-1,i\}$ and $s\in S(K)$, and let $\sv^j:=\R^jf_*\Q_p$;
this is a $\Q_p$-local system on $S_C$ \cite[Th.~10.5.1]{SW20}, and its fiber at $s$ is
$H^j_{\eet}(X_{s,C},\Q_p)$.

(i)$\Rightarrow$(iii). Let $r\geq i$. 
Since $X_s$ is proper, $\pi^{r,j}(X_s)$ is surjective, the map
$H^j_{\eet}(X_{s,C},\Q_p(r))\to \hk^{r,j}(X_s)$ is injective, and its image
is $\kker\pi^{r,j}(X_s)$ \cite[Cor.~6.15,  Rem.~6.23(i)]{CN5}. 
On the other hand, it follows from the exactness of~\eqref{gen0} {\rm (}in degree $j${\rm )}
and from the definition~\eqref{ducky1} of $F^r{\rm DR}^j(X)$,
that $\kker\pi^{r,j}(X)$ is the image of $H^j_{\proeet}(X_C,\Q_p(r))\to\hk^{r,j}(X)$.

  Now, if (i) holds, ${\cal E}^j(X)$ is acyclic. Hence, by Lemmas~\ref{simpl9}
and~\ref{simpl1}, and by~\eqref{fil}, $j_s^\dual$ induces isomorphisms
$$\hk^{r,j}(X)\simto\hk^{r,j}(X_s),\quad {\rm DR}^j(X)\simto{\rm DR}^j(X_s),\quad
F^r{\rm DR}^j(X)\simto F^r{\rm DR}^j(X_s),$$
and hence  also $\kker\pi^{r,j}(X)\simto\kker\pi^{r,j}(X_s)$.
By functoriality of the comparison maps, this implies that the map
$$j_s^\dual:H^j_{\proeet}(X_C,\Q_p(r))\lra H^j_{\eet}(X_{s,C},\Q_p(r))$$
is surjective (use \eqref{gen0}). 
But this map factors
 as
$$H^j_{\proeet}(X_C,\Q_p(r))\lra H^0(S_C,\sv^j(r))\lra H^j_{\eet}(X_{s,C},\Q_p(r)),$$
where the first map is the edge map of the Leray spectral sequence, and the second
is the evaluation at $s$. It follows that every element of the fiber of $\sv^j$ at $s$
extends to a global section of $\sv^j$ on $S_C$; since $S_C$ is connected,
this means that $\sv^j$ is constant.
(If $S$ is a dagger disc, $H^j_{\proeet}(X_C,\Q_p(r))$ is the colimit of
the $H^j_{\proeet}(X'_C,\Q_p(r))$, for $X'=\sx\times_US'$ and $S'$ running through the discs
centered at $s_0$ of radius strictly bigger than that of $S$; as $H^j_{\eet}(X_{s,C},\Q_p(r))$ is finite dimensional,
the above argument shows that $\sv^j$ is constant on $S'_C$, if $S'$ is small enough.)

(iii)$\Rightarrow$(ii).  By~\cite[Th.~8.8]{Sch13}, the local system $\R^jf_*\Q_p$ on $S$
is de Rham and the filtered module with integrable connection associated to it,
in the sense of~\cite[Def.~7.5]{Sch13}, is $({\cal H}^j,\nabla,F^\bullet)$.
Assume that $\sv^j$ is constant on $S_C$ {\rm (}replace $S$ by $S'$ in what follows,
if $S$ is a dagger disc{\rm )}. Since $S$ is geometrically connected, $\R^jf_*\Q_p$ is then
the pullback to $S$ of the representation
$V:=H^0(S_C,\sv^j)\simeq H^j_{\eet}(X_{s,C},\Q_p)$ of $G_K$.
This representation is de Rham; hence $\R^jf_*\Q_p$ is also associated to
$\O_S\otimes_KD_{\dr}(V)$, endowed with the connection $d\otimes 1$ and the
filtration $\O_S\otimes_KF^\bullet D_{\dr}(V)$.
By uniqueness of the associated filtered module with integrable
connection~\cite[Th.~7.6]{Sch13}, we have
$F^a{\cal H}^j\simeq \O_S\otimes_KF^a D_{\dr}(V)$, for all $a$.
Hence $F^a{\cal H}^j$ is generated by horizontal sections; in particular, it is
stable by $\nabla$, and $\theta_j^a=0$ for all $a$.
\end{proof}

\begin{remark}\label{etale-face}
For the family of Section~\ref{family}, Proposition~\ref{prop:etale} shows
directly that condition (iii) of Theorem~\ref{simpl2} is not satisfied for $i=1$:
the space $H^0(\So_C,\R^1f_*\Q_p)$ is of dimension~$1$ and not~$2$. This is the
\'etale counterpart of Lemma~\ref{lem:KS}.
\end{remark}

\begin{remark}[The period map]\label{periods}
Theorem~\ref{simpl2} has a simple geometric interpretation.
Let $S$ and $X$ be as in Theorem~\ref{simpl2}, and let
$V^i:=H^0(S,{\cal H}^i)^{\nabla=0}$. By Lemma~\ref{simpl9}, we have
$H^i_{\dr}(X)\simeq V^i$ and, for every $u\in S(K)$, evaluation at $u$ is an isomorphism
$V^i\simto H^i_{\dr}(X_u)$ {\rm (}horizontal transport{\rm )}. We denote by $F^\bullet_uV^i$
the inverse image by this isomorphism of the Hodge filtration of $H^i_{\dr}(X_u)$.
So all the fibers have the same de Rham cohomology $V^i$, but their Hodge filtrations
$F^\bullet_uV^i$ move with $u$; the dimension of $F^a_uV^i$ does not depend on $u$,
since $F^a{\cal H}^i$ is a subbundle of ${\cal H}^i$.
The map $u\mapsto F^\bullet_uV^i$, from $S(K)$ to the variety of flags of $V^i$
of the corresponding type, is the period map of $({\cal H}^i,\nabla,F^\bullet)$.

(i) {\it The Hodge filtration of the total space is the part common to the Hodge
filtrations of all the fibers}: for all $a$, we have
\begin{equation}\label{inter1}
F^aH^i_{\dr}(X)=\bigcap_{u\in S(K)}F^a_uV^i\quad{\text{in $V^i$.}}
\end{equation}
Indeed, by Lemma~\ref{simpl9}, $v\in V^i$ belongs to $F^aH^i_{\dr}(X)$ if and only if,
viewed as a section of ${\cal H}^i$, it is a section of $F^a{\cal H}^i$, i.e.~if and only if
its image in $H^0(S,{\cal H}^i/F^a{\cal H}^i)$ is $0$. As ${\cal H}^i/F^a{\cal H}^i$ is
locally free, and as an analytic function on $S$ which vanishes at all points of $S(K)$
is $0$, this is the case if and only if $v(u)\in F^aH^i_{\dr}(X_u)$ for all $u\in S(K)$,
i.e.~if and only if $v\in F^a_uV^i$ for all $u$. 

Since $V^i$ is finite dimensional, finitely many $u$ suffice in~\eqref{inter1}; see the
Appendix for the $\B_{\dr}^+$-version of~\eqref{inter1}.

(ii)  Consider the family of Section~\ref{family}. Let $\omega_0,\eta_0$ be the classes
of $\omega,\eta$ in $H^1_{\dr}(E)$, and identify $V^1$ with $H^1_{\dr}(E)$ by evaluation
at $T=0$. Let $G=\matrice{\alpha}{\beta}{\gamma}{\delta}$ be as in
Lemma~\ref{convergence}, viewed as a function of $\lambda$ via $U=\lambda-a$. A horizontal section $c_1(\alpha\omega+\gamma\eta)+
c_2(\beta\omega+\delta\eta)$ has value $c_1\omega_0+c_2\eta_0$ at $T=0$, and its value
at $\lambda$ belongs to $F^1=K\omega$ if and only if $c_1\gamma(\lambda)+c_2\delta(\lambda)=0$.
Hence the period map is
$$\lambda\longmapsto F^1_\lambda V^1=K\big(\delta(\lambda)\,\omega_0-\gamma(\lambda)\,\eta_0\big).$$
From $G'=-AG$ and $G(a)=I_2$, we get $\gamma(a)=0$, $\delta(a)=1$ and
$\gamma'(a)=-\frac{1}{2a(a-1)}\neq 0$ {\rm (}this is, up to sign, the coefficient of
the Kodaira--Spencer map of Lemma~\ref{lem:KS}{\rm )}. So $\gamma$ is not identically $0$, and
$F^1_\lambda V^1\neq K\omega_0=F^1_aV^1$ if $\gamma(\lambda)\neq 0$.
Two distinct lines of a plane meet in $0$, and~\eqref{inter1} gives
$F^1H^1_{\dr}(X)=0$: this is Lemma~\ref{ducky4}.
\end{remark}
\begin{corollary}\label{period-crit}
Let $S$ and $X$ be as in Hypothesis~\ref{disc}, and let $i\geq 1$. Then
Conjecture~\ref{main0} holds for $X$ in degree $i$ if and only if the period maps
of ${\cal H}^{i-1}$ and ${\cal H}^i$ are constant on $S(K)$. In particular, $X$
satisfies Conjecture~\ref{main0} if and only if all its period maps are constant.
\end{corollary}
\begin{proof}
We use the notation of Remark~\ref{periods}.
Conjecture~\ref{main0} holds for $X$ in degree $i$ if and only if ${\cal E}^j(X)$
is acyclic for $j=i-1,i$ (see Section~\ref{acyc}). Fix such a $j$, and $s\in S(K)$.
Via the isomorphisms $H^j_{\dr}(X)\simeq V^j\simeq H^j_{\dr}(X_s)$ of
Lemma~\ref{simpl9}, the map $j_s^\dual$ becomes the identity of $V^j$; hence, by
Lemma~\ref{simpl1}, ${\cal E}^j(X)$ is acyclic if and only if
$F^aH^j_{\dr}(X)=F^a_sV^j$ for all $a$. If this is the case, \eqref{inter1} gives
$F^a_sV^j\subset F^a_uV^j$, hence $F^a_uV^j=F^a_sV^j$ since both spaces have the
same dimension, for all $u\in S(K)$ and all $a$: the period map of ${\cal H}^j$ is
constant. Conversely, if it is constant, \eqref{inter1} gives
$F^aH^j_{\dr}(X)=F^a_sV^j$ for all $a$.
\end{proof}

In this language, the equivalence (i)$\Leftrightarrow$(ii) of Theorem~\ref{simpl2}
reflects the familiar fact that the Higgs field is the differential of the period map.
For the family of Section~\ref{family}, the period map of ${\cal H}^1$ is not constant
by Remark~\ref{periods}(ii), and Corollary~\ref{period-crit} shows again that $X$
contradicts Conjecture~\ref{main0}, in degrees~$1$ and~$2$.
\begin{corollary}\label{abelian}
In the situation of Theorem~\ref{simpl2}, assume that
$\sx\to U$ is an abelian scheme.  Then $X$ satisfies
Conjecture~\ref{main0} if and only if the Kodaira--Spencer map
\[
 \mathrm{KS}\colon \quad f_*\Omega^1_{X/S}\longrightarrow
 \R^1f_*\so_X\otimes_{\so_S}\Omega^1_{S/K}
\]
vanishes identically on $S$.  If $\sx\to U$ is a family of elliptic
curves, this holds if and only if the $j$-invariant is constant on $S$.
\end{corollary}

\begin{proof}
Cup product gives isomorphisms $\wedge^i\mathcal H^1\simto\mathcal H^i$
compatible with the Gauss--Manin connections and, by $E_1$-degeneration,
strictly compatible with the Hodge filtrations.  Hence, by the Leibniz
rule, $\theta_i$ is the derivation of $\wedge^i\mathcal H^1$ induced by
$\theta_1$.  So $\theta_i|_S=0$ for all $i$ if and only if
$\theta_1|_S=0$.  The only possibly nonzero component of $\theta_1$ is
$\theta^1_1\colon f_*\Omega^1_{X/S}\to\R^1f_*\so_X\otimes\Omega^1_{S/K}$,
which is $\mathrm{KS}$ (see the proof of Lemma~\ref{lem:KS}).  
Theorem~\ref{simpl2} allows us to conclude.

Suppose now that the fibers are elliptic curves. \'Etale locally on \(S\), the family has Legendre form with parameter \(\lambda\). By the formula in the proof of Lemma~\ref{lem:KS}, \(\mathrm{KS}=0\) if and only if \(d\lambda=0\), that is, if and only if \(\lambda\) is locally constant. Since the classical rational function \(j=J(\lambda)\) has finite fibers, \(\lambda\) is locally constant if and only if \(j\) is locally constant. As \(S\) is connected, this is equivalent to \(j\) being constant on \(S\).
\end{proof}

\begin{remark}
(i) The vanishing of the Kodaira--Spencer map is not a necessary condition
 for the total space of a family to satisfy
Conjecture~\ref{main0}.
For example, 
consider the whole algebraic Legendre family 
$f\colon\mathcal E\longrightarrow
 U=\mathbb P^1_{\Q_p}\setminus\{0,1,\infty\}$
of Remark~\ref{comal1.2}.
With the notations of
Lemma~\ref{lem:KS}, the same computation gives
$\mathrm{KS}(\omega)= \frac{\overline\eta}{2\lambda(\lambda-1)}\otimes d\lambda$,
which vanishes nowhere on $U$.
On the other hand, as we have seen in Remark~\ref{comal1.2},
the total space \(\mathcal E_C^{\mathrm{an}}\)
satisfies Conjecture~\ref{main0}.

(ii) The vanishing of the Higgs fields does not imply that $X\to S$ is
isotrivial.  For elliptic curves it does, by
Corollary~\ref{abelian}.  But if the fibers satisfy $h^{a,b}=0$ for
$a\neq b$, as for smooth families of del Pezzo surfaces, all Higgs maps
vanish automatically.  Every such family, isotrivial or not, then
satisfies Conjecture~\ref{main0} over the discs of
Hypothesis~\ref{disc}.
\end{remark}

\subsection{Higher dimensional bases}\label{highdim}
 Let now $U$ be a smooth rigid analytic variety of arbitrary dimension, $f:\sx\to U$
smooth and proper, $S\subset U$ a dagger or open polydisc centered at $s_0\in U(K)$,
small enough for \eqref{G-triv} to hold, and $X=\sx\times_US$. The period maps are
defined as in Remark~\ref{periods}, and the first part of the proof of Lemma~\ref{simpl9}
and Section~\ref{acyc} still give, for all $j$ and $a$, inclusions
\begin{equation}\label{highdim0}
F^aH^j_{\dr}(X)\ \subset\ H^0(S,F^a{\cal H}^j)^{\nabla=0}=\bigcap_{u\in S(K)}F^a_uV^j
\ \subset\ F^a_sV^j,
\end{equation}
and Conjecture~\ref{main0} holds for $X$ in degree $i$ if and only if the composite
inclusion is an equality for $j=i-1,i$ and all $a$. The second inclusion is an equality
for all $a$ if and only if the period map of ${\cal H}^j$ is constant, i.e., if and only if
the Higgs maps $\theta^a_j$ are $0$.  The first one is an equality if and only if all
horizontal sections of $F^a{\cal H}^j$ come from $H^j(X,\Omega^{\geq a}_X)$; this is
automatic if $S$ is a curve, by Lemma~\ref{lift}, which is why
Corollary~\ref{period-crit} holds, but it can fail if $\dim S\geq 2$.
Similarly, the implications (i)$\Rightarrow$(iii)$\Rightarrow$(ii) of
Theorem~\ref{simpl2} still hold, with the same proofs, but not their converses, as the
following example shows. 

\vskip2mm
  Let $Y$ be a smooth projective curve of genus~$2$ over $K$, with $Y(K)\neq\emptyset$,
let $L_0$ be a line bundle of degree~$\neq 0$ on $Y$, and let $q\in K^\times$ with $v_p(q)>0$.
Let $S$ be a small enough dagger polydisc centered at $0$ in the Jacobian ${\rm Pic}^0(Y)$ of $Y$ 
(in particular, $S$ is of dimension~$2$),
let $B:=Y\times S$, let $\sll$ be the restriction to $B$ of ${\rm pr}_Y^*L_0\otimes{\cal P}$,
where ${\cal P}$ is the Poincar\'e bundle on $Y\times{\rm Pic}^0(Y)$, normalized along
$\{y_0\}\times{\rm Pic}^0(Y)$ for some $y_0\in Y(K)$, 
and let $\sll^\times$ be the complement
of the zero section in the total space of $\sll$. Set
$$X:=\sll^\times/q^{\Z}\ \xrightarrow{\ \pi\ }\ B\ \xrightarrow{\ {\rm pr}_S\ }\ S,
\qquad f:={\rm pr}_S\circ\pi.$$
 Then $f$ is smooth and proper, and its fibers $X_s$ are principal bundles over $Y$ under
the Tate curve $E_q:=\mathbf G_m/q^{\Z}$. 
We claim that all the Higgs maps of $f$ are $0$, i.e., that its period maps are constant,
and that the local systems $\R^1f_*\Q_p$ and $\R^2f_*\Q_p$ are constant on $S_C$, but that
$F^2H^2_{\dr}(X)=0$ whereas $\dim_KF^2H^2_{\dr}(X_s)=2$. So, for $i=2$, conditions (ii)
and (iii) of Theorem~\ref{simpl2} and the condition of Corollary~\ref{period-crit} hold,
but the first inclusion in~\eqref{highdim0} is strict for $j=a=2$, and $X$ contradicts
Conjecture~\ref{main0} in degrees~$2$ and~$3$. 

     Let $z$ be a local coordinate on the fibers of $\sll^\times\to B$. Then
$\Omega^1_{X/B}=\O_X\frac{dz}{z}$, and the class of the extension
$0\to\pi^*\Omega^1_B\to\Omega^1_X\to\O_X\frac{dz}{z}\to 0$ is the pullback of 
$c_1(\sll)\in H^1(B,\Omega^1_B)$. Now,
$$c_1(\sll)=c\otimes 1+\tau_1\otimes ds_1+\tau_2\otimes ds_2\ \in\
\big(H^1(Y,\Omega^1_Y)\otimes_K\O(S)\big)\oplus\big(H^1(Y,\O_Y)\otimes_K\Omega^1(S)\big),$$
 where $c=c_1(L_0)\neq 0$, and where $\tau_1,\tau_2$ is a basis of $H^1(Y,\O_Y)$ and
$s_1,s_2$ are coordinates on $S$ 
{\rm (}the last two terms are the differential of the
inclusion of $S$ in ${\rm Pic}^0(Y)$, whose tangent space is $H^1(Y,\O_Y)${\rm )}.
The relative $2$-forms on $X/S$ are the $\alpha\wedge\frac{dz}{z}$, for
$\alpha\in H^0(Y,\Omega^1_Y)\otimes_K\O(S)$, and such a form lifts to a $2$-form on $X$
if and only if there exists $\beta\in\Omega^1(S)$ such that
$(\alpha+\beta)\cup c_1(\sll)=0$ in $H^1(B,\Omega^2_B)$. Writing
$\alpha\cup\tau_k=\langle\alpha,\tau_k\rangle\,c$, this means
\begin{equation}\label{highdim1}
\langle\alpha,\tau_1\rangle\,ds_1+\langle\alpha,\tau_2\rangle\,ds_2=\pm\beta
\quad{\rm and}\quad \beta\wedge ds_1=\beta\wedge ds_2=0.
\end{equation}

$\bullet$ If we restrict $f$ to a curve $S'\subset S$, the second condition is empty
{\rm (}$\Omega^2_{S'}=0${\rm )} and the first defines~$\beta$: all relative $2$-forms lift, which
means that the Higgs map $\theta_2^2$ is $0$ on $S'$. This being true for all $S'$, we
have $\theta_2^2=0$.

$\bullet$ Over $S$, the second condition implies $\beta=0$, and then the first implies
$\alpha=0$, by Serre duality on~$Y$. Hence all $2$-forms on $X$ are pullbacks of $2$-forms
on $B=Y\times S$, and their restrictions to the fibers $X_s$ are $0$ since $Y$ is a curve.
So $F^2H^2_{\dr}(X)=0$, whereas $F^2H^2_{\dr}(X_s)=H^0(X_s,\Omega^2_{X_s})\simeq
H^0(Y,\Omega^1_Y)$ is of dimension~$2$.

$\bullet$ Since $c\neq 0$, the relative $1$-forms on $X/S$ are the pullbacks of the
$\alpha\in H^0(Y,\Omega^1_Y)\otimes_K\O(S)$; they lift to $X$, hence $\theta_1^1=0$.
The two remaining Higgs maps, $\theta_2^1$ and $\theta_3^2$, are the transposes of
$\theta_2^2$ and $\theta_1^1$ for Serre duality.

$\bullet$ Since $\pi$ is a principal bundle under $E_q$, the sheaves $\R^b\pi_*\Q_p$
are the constant sheaves $H^b_{\eet}(E_{q,C},\Q_p)$; hence the $E_2$-terms of the
spectral sequence $\R^a{\rm pr}_{S*}\R^b\pi_*\Q_p\Rightarrow\R^{a+b}f_*\Q_p$ are the
constant sheaves $H^a_{\eet}(Y_C,\Q_p)\otimes H^b_{\eet}(E_{q,C},\Q_p)$, and its
differentials $d_2$ are given by cup product with $c_1(L_0)\neq 0$. One deduces that
$\R^2f_*\Q_p\simeq H^1_{\eet}(Y_C,\Q_p)\otimes H^1_{\eet}(E_{q,C},\Q_p)$, and that
$\R^1f_*\Q_p$ is an extension of $\Q_p$ by $H^1_{\eet}(Y_C,\Q_p)$, split by the class of
the $\Z$-covering $\sll^\times\to X$. So both are constant.

\vskip2mm
     In this example,  the fibers $X_s$ are $p$-adic analogues of non-K\"ahler elliptic surfaces: the same
construction with $Y=\mathbf P^1$ and $L_0=\O(-1)$ gives the Hopf surface
$(\mathbf A^2\moins\{0\})/q^{\Z}$, and with $Y$ of genus~$1$ a primary Kodaira surface.
As for these, Hodge symmetry fails: $h^{1,0}(X_s)=2$ and $h^{0,1}(X_s)=3$.
Both this failure and that of Conjecture~\ref{main0} come from
$c_1(L_0)\neq 0$, through the differential $d_2$ of the Leray spectral sequence
of~$\pi$. We do not know if the implication (ii)$\Rightarrow$(i)
of Theorem~\ref{simpl2}  can also fail when
$f$ is projective.
\appendix
\section{De Rham filtrations} We compare the Hodge filtration with the induced filtration on \(\B_{\mathrm{dR}}^+\)-cohomology and describe the associated period lattices as intersections of fiber lattices.
\subsection{The Hodge filtration}
Our goal here is to prove the following result:
\begin{proposition}\label{appendice}
If $X$ is a smooth dagger variety defined over $K$ with $H^*_{\dr}(X)$
finite dimensional, and if $n\leq r$,
$${\rm Im}\big(H^n(F^r(\rg_{\dr}(X)\otimes_K^{\LL_{\Box}}\B_{\dr}^+))\to H^n_{\dr}(X)\otimes_K\B^+_{\dr}\big)=
F^r( H^n_{\dr}(X)\otimes_K\B_{\dr}^+)$$
\end{proposition}
\begin{proof}
Let $(e_i)_{i\in I}$ be an orthonormal basis of
$\O_C$ over $\O_K$ (i.e., the images modulo a uniformizer $\pi_K$ of $K$
form an algebraic basis of
$\O_C/\pi_K\O_C$ over $\O_K/\pi_K$) and, 
if $i\in I$, choose $\hat e_i\in{\bf A}_{{\rm inf},K}:=
\O_K\otimes {\bf A}_{\rm inf}$
such that $\theta(\hat e_i)=e_i$. Then $\B_{\dr}^+$ is,
as a $K$-module\footnote{See~\cite{swan}, Lemme 3.5 and Prop.\,3.7,
for similar constructions.}, the set of
$$x=\sum_{k\geq 0}\big(\sum_{i\in I}x_{k,i}\hat e_i\big) t^k,
\quad x_{k,i}\to 0\ {\text{when $i\to \infty$}}.$$
Moreover, this decomposition of $x$ is unique, 
 $x\mapsto x_{k,i}$ is $K$-linear and continuous\footnote{We will not need
this point, but it is also true that
$x\mapsto \sum_{k\geq 0}x_{k,i}t^k$ is $K[[t]]$-linear and continuous.}, 
and $x\in F^r\B_{\dr}^+$ if and only if $x_{k,i}=0$ if $k<r$.

Now, let ${\cal U}=(U_j)$ be a covering of $X$ by dagger affinoids.  If $J$ is a finite
set of $j$'s, let $U_J:=\cap_{j\in J}U_j$.
We write an element of $Z^n({\cal U},\rg_{\dr}(X)\otimes^{\LL_{\Box}}_K\B_{\dr}^+)$ as
$(\omega_{J}^{(n_J)})_J$, with 
$n_J=n-|J|+1$ and $\omega_{J}^{(n_J)}\in\Omega^{n_J}(U_J)\otimes^{\Box}_K\B_{\dr}^+$.

The above decomposition of elements of $\B_{\dr}^+$ allows
one to decompose each $\omega_{J}^{(n_J)}$ as 
$$\omega_{J}^{(n_J)}=
\sum_{k\geq 0}\big(\sum_{i\in I}\omega_{J,k,i}^{(n_J)}\hat e_i\big) t^k$$
and $\omega_{J}^{(n_J)}\mapsto \omega_{J,k,i}^{(n_J)}$ is $K$-linear and
continuous, and commutes with differentiation and restriction to a sub-affinoid.
It follows that each $(\omega_{J,k,i}^{(n_J)})_J$ is an $n$-cocycle
belonging to $Z^n({\cal U},\rg_{\dr}(X))$.
Moreover, $(\omega_{J}^{(n_J)})_J\in Z^n(F^r(\rg_{\dr}(X)\otimes_K^{\LL_{\Box}}\B^+_{\dr}))$
if and only if $\omega_{J,k,i}^{(n_J)}=0$ if $n_J+k<r$.
Since $n_J+k<r$ is equivalent to $n_J<r-k$, this last condition
is equivalent to $(\omega_{J,k,i}^{(n_J)})_J\in Z^n({\cal U},\Omega^{\geq r-k})$
and to $t^k\hat e_i\otimes (\omega_{J,k,i}^{(n_J)})_J\in 
 Z^n({\cal U},\Omega^{\geq r-k})\otimes_K^{\LL_{\Box}}t^k\B_{\dr}^+$.
It follows that the image of $t^k\hat e_i\otimes (\omega_{J,k,i}^{(n_J)})_J$ in
$ H^n_{\dr}(X)\otimes_K\B^+_{\dr}$ belongs to $F^r(H^n_{\dr}(X)\otimes_K\B^+_{\dr})$,
which proves the inclusion $\subset$ in the proposition since $F^r( H^n_{\dr}(X)\otimes_K\B^+_{\dr})$
is closed in $ H^n_{\dr}(X)\otimes_K\B^+_{\dr}$.

In the other direction, if $x=a\otimes y$, with $a\in t^k\B_{\dr}^+$ and 
$y\in F^{r-k}H^n_{\dr}(X)$, one chooses a lifting
$(\omega_J^{(n_J)})_J$ of $y$ in $Z^n({\cal U},\Omega^{\geq r-k})$,
such that $\omega_J^{(n_J)}=0$ if $n_J<r-k$.
Then 
$$(a\omega_J^{(n_J)})_J\in Z^n({\cal U},\Omega^{\geq r-k}\otimes^{\LL_{\Box}}_Kt^k\B^+_{\dr})\subset
Z^n({\cal U},F^r(\rg_{\dr}(X)\otimes^{\LL_{\Box}}_K\B^+_{\dr})).$$
This proves the reverse inclusion, and finishes the proof.
\end{proof}
\subsection{$\B^+_{\dr}$-filtration} 
Let $S$ and $X$ be as in Theorem~\ref{simpl2}, and let
$V^i:=H^0(S,{\cal H}^i)^{\nabla=0}$, so that $H^i_{\dr}(X)\simeq V^i$ (Lemma~\ref{simpl9}).
If  $u\in S(K)$, evaluation at $u$ is an isomorphism $V^i\simto H^i_{\dr}(X_u)$.
 We have a formula similar to \eqref{inter1}  for the corresponding lattices: in $V^i\otimes_K\B_{\dr}$,
we have
\begin{equation}\label{inter2}
F^0(H^i_{\dr}(X)\otimes_K\B_{\dr})=\bigcap_{u\in S(K)}F^0(H^i_{\dr}(X_u)\otimes_K\B_{\dr}).
\end{equation}
{\rm (}If $F^\bullet$ is a decreasing, separated and exhaustive filtration of a finite
dimensional $K$-vector space $V$, then
$F^0(V\otimes_K\B_{\dr})=M(F^\bullet)\otimes_{K[t]}\B_{\dr}^+$,
where $M(F^\bullet):=\oplus_{a\in\Z}\,t^{-a}F^aV\subset V\otimes_KK[t,t^{-1}]$.
Now, $F^\bullet\mapsto M(F^\bullet)$ commutes with intersections,
and so does $M\mapsto M\otimes_{K[t]}\B_{\dr}^+$ with finite intersections, because
$\B_{\dr}^+$ is flat over $K[t]$. Hence~\eqref{inter2} follows from~\eqref{inter1}.{\rm )}
 
    In other words, if we view the ${\cal E}^i(X_u)$ and ${\cal E}^i(X)$ as modifications
at $x_\infty$ of the vector bundle attached to $H^i_{\hk}(X_C)$, we have
$${\cal E}^i(X)=\bigcap_{u\in S(K)}{\cal E}^i(X_u),\qquad
\deg{\cal E}^i(X)-\deg{\cal E}^i(X_s)
=-\sum_{a\geq 1}\big(\dim_KF^aH^i_{\dr}(X_s)-\dim_KF^aH^i_{\dr}(X)\big).$$
Here $s\in S(K)$ is arbitrary (the dimensions
$\dim_K F^a H^i_{\dr}(X_s)$ are independent of $s$). 
Each of the ${\cal E}^i(X_u)$ is pure of slope~$0$ (hence also of degree~$0$), 
and their intersection is of slope~$0$
if and only if they are all equal: this is another way to phrase
Lemma~\ref{simpl1}. Also, if $r\geq i$, we obtain, inside $\hk^{r,i}(X)$,
$$\kker\pi^{r,i}(X)=\bigcap_{u\in S(K)}\kker\pi^{r,i}(X_u)=
\bigcap_{u\in S(K)}H^i_{\eet}(X_{u,C},\Q_p(r)).$$

\end{document}